\documentclass[11pt]{article}
\usepackage{amssymb}
\usepackage{amsmath}
\usepackage{amsbsy}
\usepackage{epsfig}
\usepackage{enumerate}
\usepackage{bm}

\newtheorem{theorem}{Theorem}[section]

\newtheorem{condition}{Condition}

\newtheorem{corollary}{Corollary}[section]
\newtheorem{remark}{\sc Remark}[section]

\renewcommand{\theequation}{\thesection.\arabic{equation}}

\def\proclaim#1{\par \bigskip\noindent {\bf #1}\bgroup\it\ }
\def\endproclaim{\egroup\par\bigskip}

\newbox\TempBox \newbox\TempBoxA

\newcommand{\non}{\nonumber \\} 
\def\pr{\textsf{P}} 
\def\ep{\textsf{E}} 
\def\Var{\textsf{Var}} 
\def\logit{\textrm{logit}}
\def\Cal#1{{\mathcal #1}}
\def\bk#1{\bm #1}
\def\bkg#1{\bm #1} 
\def\smallbkg#1{\bm #1} 
\def\underwiggle 1{
\ifmmode\setbox\TempBox=\hbox{$ 1$}\else\setbox\TempBox=\hbox{1}\fi
\setbox\TempBoxA=\hbox to \wd\TempBox{\hss\char'176\hss}
\rlap{\copy\TempBox}\smash{\lower9pt\hbox{\copy\TempBoxA}}}
\begin{document}

\begin{center}
{ \large \bf ASYMPTOTIC PROPERTIES OF COVARIATE-ADJUSTED ADAPTIVE
DESIGNS}
\end{center}
\smallskip

\begin{center} {\sc \small By
Li-Xin Zhang\footnote{Research supported by a grant from the
National Natural Science Foundation of China (No. 10471126).},
Feifang Hu$\footnote{Research supported by grants DMS-0204232 and
DMS-0349048 from the National Science Foundation (USA).}$, Siu Hung
Cheung$^3$ and Wai Sum Chan\footnote{ Research supported by a grant
from the Research Grants Council of the Hong Kong Special
Administrative Region (Project no. CUHK400204).).
\\ \hspace*{0.33cm}
\\ \hspace*{0.5cm}
 }}
\bigskip

{\sl \small Zhejiang University, University of Virginia, Chinese
University of Hong Kong and University of Hong Kong  }
\end{center}

\begin{abstract}
Response-adaptive designs have been extensively studied
and used in clinical trials.
However, there is a lack of a comprehensive study of
response-adaptive designs that include covariates,
 despite their importance in clinical experiments.
Because the allocation scheme and the estimation of parameters are affected
by both the responses and the covariates, covariate-adjusted response-adaptive (CARA)
designs are very complex to formulate.
In this paper, we overcome
the technical hurdles and
lay out a framework for general CARA designs
for the allocation of subjects to $K (\geq 2)$ treatments.
The asymptotic properties
are studied under certain widely satisfied conditions. The
proposed CARA designs can be applied
to generalized linear models.  Two important
special cases, the linear model and the logistic regression model, are
considered in detail.
\end{abstract}

\baselineskip 14.9pt
\setcounter{section}{1}
\setcounter{equation}{0}
\setcounter{theorem}{0}
\setcounter{remark}{0}

{\bf 1. Preliminaries.~~}

1.1. {\em Brief history.}~~ In most clinical trials, patients
accrue sequentially.  Response-adaptive designs provide allocation
schemes that assign different treatments to incoming patients
based on the previous observed responses of patients. A major objective
of response-adaptive designs in clinical trials is to minimize
the number of patients that is assigned to the inferior treatment to a degree that still
generates useful statistical inferences. The ethical and other characteristics of
response-adaptive designs have been extensively discussed by many
authors ({\em e.g.}, Zelen and Wei (1995)).

Early important work on response-adaptive designs was carried out by
 Thompson (1933) and Robbins (1952).  Since then, a
steady stream of research (Zelen (1969), Wei and Durham
(1978), Eisele and Woodroofe (1995)) in this area has
generated various treatment allocation schemes for clinical
trials. Some of the advantages of using response-adaptive designs
have been recently studied by Hu and Rosenberger (2003) and
Rosenberger and Hu (2004).

In many clinical trials, covariate information is available that
has a strong influence on the responses of patients. For instance, the
efficacy of a hypertensive drug is related to a patient's initial
blood pressure and cholesterol level, whereas the effectiveness of a
cancer treatment may depend on whether the patient is a smoker or
a non-smoker.

The following notations and definitions are introduced to describe
the randomized treatment allocation schemes. Given a clinical
trial with $K$ treatments. Let $\bk X_1, \bk X_2,...$ be the
sequence of random treatment assignments.  For the $m$-th subject,
 $\bk
X_m=(X_{m,1},\ldots,X_{m,K})$ represents the assignment of
treatment such that if the $m$-th subject is allocated to
treatment $k$, then all elements in $\bk X_m$ are $0$ except for
the $k$-th component, $X_{m,k}$, which is $1$. Let $N_{m,k}$ be
the number of subjects assigned to treatment $k$ in the first $m$
assignments and write $\bk N_m=(N_{m,1},\ldots, N_{m,K})$. Then
$\bk N_m=\sum_{i=1}^m \bk X_i$. Suppose that $\{ Y_{m,k},
~k=1,\ldots, K,~ m=1,2\ldots\}$ denote the responses such that
$Y_{m,k}$ is the response of the $m$-th subject to treatment $k$,
$k=1,\ldots, K$.  In practical situations, only $Y_{m,k}$ with
$X_{m,k}=1$ is observed. Denote $\bk Y_m=(Y_{m,1},\ldots,
Y_{m,K})$. Let $\Cal X_m=\sigma(\bk X_1,\ldots,\bk X_m)$ and $\Cal
Y_m=\sigma(\bk Y_1,\ldots,\bk Y_m)$ be the corresponding sigma
fields. A response-adaptive design is defined by
$$\bkg \psi_m=E(\bk X_m|\Cal X_{m-1}, \Cal Y_{m-1}).$$

Now, assume that covariate information is available in the
clinical study. Let $\bkg\xi_m$ be the covariate of the $m$-th
subject and $\Cal Z_m=\sigma(\bkg\xi_1,\ldots,\bkg\xi_m)$ be the
corresponding sigma field. In addition, let $\Cal F_m=\sigma(\Cal
X_m, \Cal Y_m, \Cal Z_m)$ be the sigma field of the history. A
general covariate-adjusted response-adaptive (CARA) design is
defined by
$$\bkg \psi_m=E(\bk X_m|{\cal F}_{m-1}, \bkg\xi_m)=E(\bk X_m| \Cal X_{m-1}, \Cal Y_{m-1}, \Cal Z_m),$$
the conditional probabilities of assigning treatments $1,...,K$ to
the $m$th patient, conditioning on the entire history including
the information of all previous $m-1$ assignments, responses, and
covariate vectors, plus the information of the current patient's
covariate vector.

A number of attempts have been made to formulate adaptive designs
in the presence of covariates. For example, Zelen (1974) and
Pocock and Simon (1975) considered balancing covariates by using
the idea of the biased coin design (Efron, 1971). Atkinson (1982,
1999, 2002) tackled this problem by employing the $D$-optimality
criterion with a linear model. The prime concern of these works is
to balance allocations over the covariates with treatment
assignment probabilities
$$\bkg \psi_m=E(\bk X_m|\Cal X_{m-1}, \Cal Z_m)$$
which differs from the CARA designs.  These allocation schemes do
not depend on the outcome of the treatment which is important for
adaptive designs that aim to reduce the number of patients that
receive the inferior treatment.

The history to incorporate covariates in response-adaptive designs
is short. For the randomized play-the-winner rule, Bandyopadhyay and
Biswas (1999) incorporated ploytomous covariates with binary
responses. Rosenberger, Vidyashankar and Agarwal (2001) considered a
CARA design for binary responses that use a logistic regression
model. Their encouraging simulation study indicates that their
approach, together with the inclusion of the covariates,
significantly reduces the percentage of treatment failures. However,
theoretical justifications and asymptotic properties have not been
given. Further, the applications of their procedure are limited to
two treatments with binary responses.

To compare two treatments, Bandyopadhyay and Biswas (2001)
considered a linear model to utilize covariate information with
continuous responses. A limiting allocation proportion was also
derived in their design. However, according to their proposed
scheme, the conditional assignment probabilities are
$$\bkg \psi_m=E(\bk X_m|{\cal F}_{m-1}).$$
The above probabilities do not incorporate the covariates of the
incoming patient, which in some cases are crucial. For instance,
let the covariate be gender and there are two treatments, and we
assume that male and female patients react very differently to
treatments A and B. Whether the next patient is male or female
should therefore be considered as an important element that
affects the assignment of treatment. Recently, Atkinson (2004)
considered adaptive biased-coin designs for $K$-treatment based on
a linear regression model. Atkinson and Biswas (2005a and 2005b)
proposed adaptive biased-coin designs and Bayesian adaptive
biased-coin designs for clinical trials with normal responses.
However, none of these articles provided asymptotic distribution
of the estimators and allocation proportion. Without the asymptotic
properties of the estimators, it is difficult to assess the
validity of the statistical inferences after using CARA designs.

Instead of working on specific setups, we seek to derive a general
framework of CARA designs and provide theoretical foundation for
using CARA design. In a CARA design, the assignment of treatment
$\bk X_m$ depends on ${\cal F}_{m-1}$ and the covariate
information ($\bkg\xi_m$)  of the incoming patient.  This
generates a certain level of technical complexity. However, it is
important to provide a solid foundation (including asymptotic
normality) for CARA designs that can be usefully applied in many
circumstances.

1.2. {\em Main results and organization of the paper.}~~
The main objectives
are (i) to propose a general CARA
design that can be applied to cases in which
$K$-treatments ($K \geq 2$) are present and to different types of
responses (discrete or continuous), and (ii) to study important
asymptotic properties of the CARA
design. These properties provide a solid foundation for both the CARA design and
the statistical inference after using a CARA design.
 Major mathematical techniques, including martingale theory
and Gaussian approximation, are employed to develop the asymptotic results.

The rest of the paper is organized as follows. In Section 2, we introduce the
general framework of the CARA design.
Useful asymptotic results (including the strong consistency and
asymptotic normality) of  both the estimators of the unknown parameters and
the allocation proportions are derived.
The generalized linear model represents a
broad class of applications and is an important tool in the analysis of
data that involve covariates.
In Section 3, the CARA
design is applied to generalized linear models, and two
important special cases, the linear model and the logistic
regression model, are considered in detail. Under the general framework,
we are able to propose many new and useful CARA designs.
We then conclude our paper with some
observations in Section 4. Technical proofs are provided in the
Appendix.

\smallskip
\setcounter{section}{2} \setcounter{equation}{0}
\setcounter{theorem}{0} \setcounter{remark}{0}

{\bf 2. General CARA design.}~~

2.1. {\em General framework.}~~
Based on the notation in Section 1, supposing that a patient with a
covariate vector $\bkg\xi$ is assigned to treatment $k$,
$k=1,\ldots, K$, and the observed response is $Y_k$. Assume that the
responses and the covariate vector satisfy
$$ \ep[Y_k|\bkg\xi]=p_k(\bkg\theta_k,\bkg\xi),~~ \bkg\theta_k\in
\bkg\Theta_k, \quad k=1,\ldots, K,
$$
where $p_k(\cdot,\cdot)$, $k=1,\ldots, K$, are known functions.
Further, $\bkg\theta_k$, $k=1,\ldots, K$, are unknown parameters,
and $\bkg\Theta_k\subset \mathbb R^d$ is the parameter space of
$\bkg\theta_k$.  Write $\bkg\theta=(\bkg\theta_1,\ldots,
\bkg\theta_K)$ and $\bkg\Theta=\bkg\Theta_1\times
\cdots\times\bkg\Theta_K$. This model is quite general, and includes the
 generalized linear models of McCullagh and Nelder
(1989) as special cases. The discussion of the generalized linear
models is undertaken in Section 3.
We assume that $\{(Y_{m,1},\ldots, Y_{m,K},\bkg\xi_m),
~m=1,2,\ldots\}$ is a sequence of i.i.d. random vectors, the
distributions of which are the same as that of $(Y_1,\ldots, Y_K,\bkg\xi)$.

2.2. {\em CARA design.}~~ The allocation scheme is as follows.  To
start, assign $m_0$ subjects to each treatment by using a restricted
randomization. Assume that $m$ ($m\ge K m_0$) subjects have been
assigned to treatments. Their responses $\{\bk Y_j, ~j=1,\ldots,m\}$
and the corresponding covariates $\{\bkg\xi_j, ~j=1,\ldots, m\}$ are
observed. We let
$\widehat{\bkg\theta}_m=(\widehat{\bkg\theta}_{m,1},\ldots,\widehat{\bkg\theta}_{m,K})$
be an estimate of $\bkg\theta=(\bkg\theta_1,\ldots,\bkg\theta_K)$.
Here, for each $k=1,\ldots, K$, $\widehat{\bkg\theta}_{m,k}
=\widehat{\bkg\theta}_{m,k}(Y_{j,k},\bkg\xi_j: X_{j,k}=1,
j=1,\ldots,m)$ is the estimator of $\bkg\theta_k$ that is based on
the observed $N_{m,k}$-size sample $\{(Y_{j,k},\bkg\xi_j):$ for
which $X_{j,k}=1, j=1\ldots,m\}$. Next, when the $(m+1)$-th subject
is ready for randomization and the corresponding covariate
$\bkg\xi_{m+1}$ is recorded, we assign the patient to treatment $k$
with a probability of
\begin{eqnarray}\label{eqassigment}
\psi_k=\pr\big( X_{m+1,k}=1 \big|\Cal F_m, \bkg\xi_{m+1}\big)
=\pi_k( \widehat{\bkg\theta}_m,\bkg\xi_{m+1})\quad k=1,\ldots,K,
\end{eqnarray}
where $\Cal F_m=\sigma(\bk X_1,\ldots,\bk X_m,\bk Y_1,\ldots,\bk
Y_m, \bkg\xi_1,\ldots,\bkg\xi_m)$ is the sigma field of the
history and $\pi_k(\cdot,\cdot)$, $k=1,...,K$ are some given
functions.  Given $\Cal F_m$ and $\bkg\xi_{m+1}$, the response
$\bk Y_{m+1}$ of the $(m+1)$-th subject is assumed to be
independent of its assignment $\bk X_{m+1}$.  We call the function
$\bkg\pi(\cdot,\cdot)=\big(\pi_1(\cdot,\cdot),\ldots,
\pi_K(\cdot,\cdot)\big)$ the allocation function that satisfies
$\pi_1+\cdots+\pi_K\equiv 1$. Let
$g_k(\bkg\theta^{\ast})=\ep[\pi_k(\bkg\theta^{\ast},\bkg\xi)]$.
From (\ref{eqassigment}), it follows that
\begin{eqnarray}
\pr\big( X_{m+1,k}=1 \big|\Cal F_m\big)
=g_k(\widehat{\bkg\theta}_m), \quad k=1,\ldots,K.
\end{eqnarray}

Different choices of $\bkg\pi(\cdot,\cdot)$ generate different
possible classes of useful designs.   For example, we can take
$\pi_k(\bkg\theta,\bkg\xi)=R_k(\bkg\theta_1\bkg\xi^{T},$ $
\ldots,\bkg\theta_K\bkg\xi^{T})$, $k=1,\ldots,K$,  which includes
a large class of applications. Here, $0<R_k(\bk z)<1 $,
$k=1,\ldots,K$, are real functions that are defined in $\mathbb
R^K$ with
\begin{eqnarray}
\label{eqfunctionR}  \quad \sum_{k=1}^K R_k(\bk
z)=1 \quad \text{and } \; R_i(\bk z)=R_j(\bk z) \; \text{
whenever }\; z_i=z_j.
\end{eqnarray}
For simplicity, it is assumed that $\bkg\xi$ and $\bkg\theta_k$,
$k=1,\ldots, K$ have the same dimensions, otherwise, slight
modifications are necessary (see Example 3.1 for an illustration).
 In practice, the functions $R_k$ can be defined as
$$ R_k(\bk z)=\frac{G(z_k)}{G(z_1)+\cdots+G(z_K)}, \quad k=1,\ldots,K, $$
where $G$ is a smooth real function that is defined in $\mathbb R$
and satisfies $0<G(z)<\infty$. An example is that $R_k(\bk z)=e^{T
z_k}/(e^{T z_1}+ \cdots+e^{T z_K})$, $k=1,\ldots,K$.

In the two-treatment case, we can let $R_1(z_1,z_2)=G(z_1-z_2)$
and $R_2(z_1,z_2)=G(z_2-z_1)$, where $G$ is real function defined
on $\mathbb R$ satisfying $G(0)=1/2$, $G(-z)=1-G(z)$ and
$0<G(z)<1$ for all $z$. For the logistic regression model,
Rosenberger, Vidyashankar and Agarwal (2001) suggested using the
estimated covariate-adjusted odds ratio to allocate subjects,
which is equivalent to defining
$R_k(z_1,z_2)=e^{z_k}/(e^{z_1}+e^{z_2})$, $k=1,2$. For each fixed
covariate $\bkg\xi$, we can also choose
$\bkg\pi(\bkg\theta,\bkg\xi)$ according to Baldi Antognini and
Giovagnoli (2004) and Hu and Rosenberger (2003). When
$\bkg\pi(\bkg\theta,\bkg\xi)$ does not depend on $\bkg\xi$, one
can use the allocation scheme of Bandyopadhyay and Biswas (2001)
for the normal linear regression model. We now introduce some
important asymptotic properties.

\smallskip
2.3. {\it Asymptotic properties.}~~ Write
$\bkg\pi(\bkg\theta^{\ast},\bk x
)=\big(\pi_1(\bkg\theta^{\ast},\bk
x),\ldots,\pi_K(\bkg\theta^{\ast},\bk x)\big)$, $\bk
g(\bkg\theta^{\ast})=\big(g_1(\bkg\theta^{\ast}),\ldots,g_K(\bkg\theta^{\ast})\big)$,
$v_k=g_k(\bkg\theta)=\ep [\pi_k(\bkg\theta,\bkg\xi)]$, $
k=1,\ldots, K$, and $\bk v=(v_1,\ldots, v_K)$. We assume that
$0<v_k<1$, $k=1,\ldots, K$. For the allocation function $ \bkg\pi
(\bkg\theta^{\ast},\bk x)$ we assume the following condition.

\begin{condition}\label{conditionA}
We assume that  the parameter space $\bkg\Theta_k$ is a bounded
domain in $\mathbb R^d$, and that the true value $\bkg\theta_k$ is an
interior point of $\bkg\Theta_k$, $k=1,\ldots, K$.
\begin{enumerate}
\item
For each fixed $\bk x$, $\pi_k(\bkg\theta^{\ast},\bk x)>0$
is a continuous function of $\bkg\theta^{\ast}$, $k=1,\ldots,K$.
\item
For each $k=1,\ldots, K$, $\pi_k(\bkg\theta^{\ast},\bkg\xi)$
is differentiable with respect to $\bkg\theta^{\ast}$ under the
expectation, and there is a $\delta>0$ such that
$$g_k(\bkg\theta^{\ast})=g_k(\bkg\theta)+(\bkg\theta^{\ast}-\bkg\theta)\big(\frac{\partial
g_k}{\partial \bkg\theta^{\ast}}\Big|_{\smallbkg\theta}\big)^{T}
+o(\|\bkg\theta^{\ast}-\bkg\theta\|^{1+\delta}),
$$
where $\partial g_k/\partial\bkg\theta^{\ast}=(\partial
g/\partial\theta^{\ast}_{11},\ldots,\partial
g/\partial\theta^{\ast}_{Kd})$.
\end{enumerate}
\end{condition}

\begin{theorem}\label{th2} Suppose that for $k=1,\ldots, K$,
\begin{eqnarray}\label{rateconsi}
\quad \widehat{\bkg\theta}_{nk}-\bkg\theta_k=\frac{1}{n}\sum_{m=1}^n
X_{m,k} \bk h_k(Y_{m,k},\bkg\xi_m)\big(1+o(1)\big)+o(n^{-1/2})
\quad a.s.,
\end{eqnarray}
where $\bk h_k$ are $K$
functions with $\ep[\bk h_k(Y_k,\bkg\xi)|\bkg\xi]=\bk 0$.
 We also
assume that $\ep\|\bk h_k(Y_k,\bkg\xi)\|^2<\infty$, $k=1,\ldots, K$.
Then under Condition \ref{conditionA}, we have for $k=1,\ldots,K$,
{\em
\begin{eqnarray}\label{eqConsiAssProb} \pr\big(X_{n,k}=1 \big)\to v_k; \quad
\pr\big(X_{n,k}=1|\Cal F_{n-1}, \bkg\xi_n=\bk x \big)\to \pi_k(\bkg\theta,
\bk x) \; a.s.
\end{eqnarray}}
and
\begin{eqnarray}\label{eqLIL}
 \frac{\bk N_n}{n}-\bk
v=O\Big(\sqrt{\frac{\log\log n}{n}}\Big) \; a.s.;\quad
 \widehat{\bkg\theta}_n-\bkg\theta=O\Big(\sqrt{\frac{\log\log
 n}{n}}\Big).
 \end{eqnarray}
Further, let
$ \bk V_k=\ep\{\pi_k(\bkg\theta,\bkg\xi)(\bk h_k(Y_k,\bkg\xi))^{T}
\bk h_k(Y_k,\bkg\xi)\}, ~ k=1,\ldots, K,$

\noindent
$ \bk V=diag\big(\bk V_1, \ldots, \bk
V_K\big),
$
$\bkg\Sigma_1=diag(\bk v)-\bk v^{T}\bk v, \; \bkg\Sigma_2=\sum_{k=1}^K \frac{\partial
\bk g}{\partial\bkg\theta_k} \bk V_k \big(\frac{\partial \bk
g}{\partial\bkg\theta_k}\big)^{T}
$
and $\bkg\Sigma=\bkg\Sigma_1+2\bkg\Sigma_2.$  Then,
\begin{eqnarray}\label{eqCLT} \sqrt{n}(\bk N_n / n-\bk v) \overset{D}\to N(\bk 0,\bkg\Sigma)
\; \text{ and } \; \sqrt{n}(\widehat{\bkg\theta}_n-\bkg\theta)
\overset{D}\to N(\bk 0,\bk V).
\end{eqnarray}
\end{theorem}

\begin{remark} Condition (\ref{rateconsi})
depends on different estimation methods. In the next  section, we
show that it is satisfied in many cases.
\end{remark}

Theorem \ref{th2} provides general results on the asymptotic
properties of the allocation proportions $N_{n,k}/n$, $k=1,\ldots,
K$. Sometimes, one may be interested in the proportions for a
given covariate (for discrete $\bkg\xi$) as discussed in Section
3. Given a covariate $\bk x$, the proportion of subjects that is
assigned to treatment $k$ is
$$ \frac{\sum_{m=1}^n X_{m,k} I\{\bkg\xi_m=\bk x\}}
{\sum_{m=1}^nI\{\bkg\xi_m=\bk x\}}:=\frac{ N_{n,k|\bk x} }{ N_n (\bk x) },
$$
where $N_{n,k|\bk x}$ is the number of subjects with covariate
$\bk x$ that is randomized to treatment $k$, $k=1,\ldots, K$, in
the $n$ trials, and $N_n(\bk x)$ is the total number of subjects
with covariate $\bk x$. Write $\bk N_{n|\bk x}=(N_{n,1|\bk
x},\ldots,N_{n,K|\bk x})$.
The following theorem establishes the  asymptotic results of these
proportions.

\begin{theorem} \label{givenCov} Given a covariate $\bk x$, suppose that
{\em $\pr\big(\bkg\xi=\bk x \big)>0$}. Under Condition \ref{conditionA} and
(\ref{rateconsi}), we have
\begin{eqnarray}\label{givenCovConsi}
 N_{n,k|\bk x}/N_n(\bk x)\to \pi_k(\bkg\theta,\bk x)\;\;
a.s. \quad k=1,\ldots,K
\end{eqnarray}
and
\begin{eqnarray}\label{givenCovNormal}
\sqrt{N_n(\bk x)}\left(\bk N_{n|\bk x}/N_n(\bk
x)-\bkg\pi(\bkg\theta,\bk x) \right)\overset{\Cal D}\to N(\bk 0,
\bkg\Sigma_{|\bk x}),
\end{eqnarray}
 where
$$\bkg\Sigma_{|\bk x}=diag(\bkg\pi(\bkg\theta,\bk x))-
\bkg\pi(\bkg\theta,\bk x)^{T}\bkg\pi(\bkg\theta,\bk x)+2\sum_{k=1}^K \frac{\partial
\bkg\pi(\bkg\theta,\bk x)}{\partial \bkg\theta_k} \bk V_k
\Big(\frac{\partial \bkg\pi(\bkg\theta,\bk x)}{\partial
\bkg\theta_k}\Big)^{T}\pr(\bkg\xi=\bk x). $$
\end{theorem}

\setcounter{section}{3} \setcounter{equation}{0}
\setcounter{theorem}{0} \setcounter{remark}{0}

{\bf 3. Generalized linear models.}~~

In this section, the general results of Section 2 are applied to
the generalized linear model (GLM) and its two special cases, the
logistic regression model and the linear model (refer to McCullagh
and Nelder (1989) for applications of these models). Suppose,
given $\bkg\xi$, that the response $Y_k$ of a trial of treatment
$k$ has a distribution in the exponential family, and takes the
form
\begin{eqnarray}\label{eqglm}
f_k(y_k|\bkg\xi,
\bkg\theta_k)=\exp\big\{(y_k\mu_k-a_k(\mu_k))/\phi_k+b_k(y_k,\phi_k)\}
\end{eqnarray}
 with link function $\mu_k=h_k(\bkg\xi\bkg\theta_k^{T})$,
 where $\bkg\theta_k=(\theta_{k1},\ldots,\theta_{kd})$,
$k=1,\ldots, K$, are coefficients. Assume that the scale parameter
$\phi_k$ is fixed, then
$ \ep[Y_k|\bkg\xi]=a_k^{\prime}(\mu_k),~
\Var(Y_k|\bkg\xi)=a_k^{\prime\prime}(\mu_k)\phi_k $ and
$$ \frac{\partial\log f_k(y_k|\bkg\xi,\bkg\theta_k)}{\partial
\bkg\theta_k}=\frac{1}{\phi_k}\{y_k-a_k^{\prime}(\mu_k)\}
h_k^{\prime}(\bkg\xi\bkg\theta_k^T)\bkg\xi,$$
$$ \frac{\partial^2\log f_k(y_k|\bkg\xi,\bkg\theta_k)}{\partial
\bkg\theta_k^2}
=\frac{1}{\phi_k}\Big\{-a_k^{\prime\prime}(\mu_k)[h_k^{\prime}(\bkg\xi\bkg\theta_k^T)]^2
+[y_k-a_k^{\prime}(\mu_k)]h_k^{\prime\prime}
(\bkg\xi\bkg\theta_k^T)\Big\}\bkg\xi^{T}\bkg\xi.$$
Thus, given $\bkg\xi$, the conditional Fisher information matrix is
$$ \bk I_k(\bkg\theta_k|\bkg\xi)=-\ep\Big[\frac{\partial^2\log
f_k(Y_k|\bkg\xi,\bkg\theta_k)}{\partial
\bkg\theta_k^2}\Big|\bkg\xi\Big]
=\frac{1}{\phi_k}a_k^{\prime\prime}(\mu_k)[h_k^{\prime}
(\bkg\xi\bkg\theta_k^T)]^2\bkg\xi^{T}\bkg\xi.$$
 For the observations up to stage
$m$, the likelihood function is
$$
L(\bkg\theta)=\prod_{j=1}^m\prod_{k=1}^K[f_k(Y_{j,k}|\bkg\xi_j,
\bkg\theta_k)]^{X_{j,k}}=\prod_{k=1}^K
\prod_{j=1}^m[f_k(Y_{j,k}|\bkg\xi_j,\bkg\theta_k)]^{X_{j,k}}:=
\prod_{k=1}^KL_k(\bkg\theta_k)$$ with $ \log
L_k(\bkg\theta_k)\propto \sum_{j=1}^m
X_{j,k}\{Y_{j,k}-a_k(\mu_{j,k})\},\;\;
\mu_{j,k}=h_k(\bkg\theta_k^{T}\bkg\xi_j), \;\; k=1,2,\ldots,K. $ The
MLE
$\widehat{\bkg\theta}_m=(\widehat{\bkg\theta}_{m,1},\ldots,\widehat{\bkg\theta}_{m,K})$
of $\bkg\theta=(\theta_1,\ldots,\theta_K)$ is that for which
$\widehat{\bkg\theta}_m$  maximizes  $L(\bkg\theta)$ over
$\bkg\theta\in \bkg\Theta_1\times\cdots\times\bkg\Theta_K$.
Equivalently, $ \widehat{\bkg\theta}_{m,k}$ maximizes $L_k$ over
$\bkg\theta_k\in\bkg\Theta_k$, $k=1,2,\ldots, K$. Rosenberger,
Flournoy and Durham (1997) established a general result for the
asymptotic normality of MLEs from a response-driven design.
Rosenberger and Hu (2002) gave the asymptotic normality of the
regression parameters from a generalized linear model that followed
a sequential design with covariate vectors.  These two papers
neither examined the case of using covariates to adjust the design,
nor established the asymptotic properties of the allocation
proportions. The next corollary gives results on both the estimators
of the parameters and the allocation proportions.
\begin{corollary}  \label{thglm}
 Define
\begin{eqnarray}\label{eqIk}\bk I_k=\bk I_k(\bkg\theta)
=\ep\{\pi_k(\bkg\theta,\bkg\xi)\bk I_k(\bkg\theta_k|\bkg\xi)\},
\;\; k=1,\ldots, K.
\end{eqnarray}
Under Condition \ref{conditionA},
if the matrices $\bk I_k$, $k=1,2\ldots, K$, are
nonsingular and the MLE $\widehat{\bkg\theta}_m$ is unique, then under
 regularity condition (A.13) in the Appendix, we have (\ref{eqConsiAssProb}),
(\ref{eqLIL}), and (\ref{eqCLT}) with $\bk V_k=\bk I_k^{-1}$,
$k=1,\ldots, K$.
Moreover, if $\pr(\bkg\xi=\bk x)>0$ for a given covariate $\bk x$,
then (\ref{givenCovConsi}) and (\ref{givenCovNormal}) hold.
\end{corollary}
This result is a corollary of Theorems \ref{th2} and
\ref{givenCov}. The proof is given  in the Appendix through
the verification of Condition (\ref{rateconsi}). For both the logistic regression and
the linear regression, condition (A.13) is satisfied.

\begin{remark}\label{rek3.1}
 From Theorem \ref{thglm}, it follows that
\begin{eqnarray}\label{eqrek3.1.1}
 \quad \quad \sqrt{N_{n,k}}\;
(\widehat{\bkg\theta}_{n,k}-\bkg\theta_k)\overset{D}\to N\big(\bk
0,v_k \{\ep[\pi_k(\bkg\theta,\bkg\xi)\bk
I_k(\bkg\theta_k|\bkg\xi)]\}^{-1}\big),\quad k=1,\ldots, K.
\end{eqnarray}
It should be noted that the asymptotic variances are different from
those of general linear models with a fixed allocation procedure. For
the latter, we have
 \begin{eqnarray}\label{eqrek3.1.2}
  \sqrt{N_{n,k}}\;
(\widehat{\bkg\theta}_{n,k}-\bkg\theta_k)\overset{D}\to
N\big(0,\{\ep[\bk I_k(\bkg\theta_k|\bkg\xi)]\}^{-1} \big), \quad
k=1,\ldots, K.
\end{eqnarray}
If the allocation functions $\pi_k(\bkg\theta,\bkg\xi)$ do not
depend on $\bkg\xi$, then
$\pi_k(\bkg\theta,\bkg\xi)=g_k(\bkg\theta)=v_k$, and so
(\ref{eqrek3.1.1}) and (\ref{eqrek3.1.2}) are identical. Our
asymptotic variance-covariance matrix of $\widehat{\bkg\theta}_n$
is also different from that in Theorem 2 of Baldi  Antognini and
Giovagnoli (2004), because the allocation probabilities in their
study  do not depend on the covariates.
\end{remark}

\begin{remark}\label{rkvarEstGLM} When the distribution of
$\bkg\xi$ and the true value of $\bkg\theta$ are known, the values
of $\bk v=\ep[\bkg\pi(\bkg\theta,\bkg\xi)]$, $\partial \bk
g/\partial \bkg\theta_k=\ep[ \partial
\bkg\pi(\bkg\theta,\bkg\xi)/\partial \bkg\theta_k]$ and $\bk I_k$
in (\ref{eqIk}) can be obtained by computing the expectations, and
then the values of the asymptotic variance-covariance matrices
$\bk V$, $\bkg\Sigma$ and $\bkg\Sigma_{|\bk x}$ can be obtained.
In practice, we can obtain the estimates as follows.

\noindent
\begin{enumerate}
\item[(a)] Estimate $\bk I_k$ by $ \widehat{\bk I}_{n,k}=\frac 1n
\sum_{m=1}^n X_{m,k}\bk
I_k(\widehat{\bkg\theta}_{n,k}|\bkg\xi_m),\quad k=1,2,\ldots, K; $
and then the estimator of $\bk V$ is $ \widehat{\bk V}_n=diag\big(
\widehat{\bk I}_{n,1}^{-1},\ldots,\widehat{\bk I}_{n,K}^{-1}\big).
$ \item[(b)] Estimate  $\bkg\Sigma_1$ and $\frac{\partial\bk
g}{\partial \smallbkg\theta_k}$, respectively, by
$$\widehat{\bkg\Sigma_1}=diag(\frac{\bk N_n}{n})-
\big(\frac{\bk N_n}{n}\big)^{T}\frac{\bk N_n}{n} \;\;\text{ and
}\;\; \widehat{\frac{\partial \bk g}{\partial
\bkg\theta_k}}=\frac{1}{n}\sum_{m=1}^n \frac{\partial
\bkg\pi(\bkg\theta^{\ast},\bkg\xi_m)}{\partial
\bkg\theta^{\ast}_k}\big|_{\smallbkg
\theta^{\ast}=\widehat{\smallbkg\theta}_n}.
$$
\item[(c)] Define the estimator $\widehat{\bkg\Sigma}$ of
$\bkg\Sigma$ by
$\widehat{\bkg\Sigma}=\widehat{\bkg\Sigma_1}+2\sum_{k=1}^K\widehat{\frac{\partial
\bk g}{\partial \bkg\theta_k}}\widehat{\bk
V}_{n,k}\big(\widehat{\frac{\partial \bk g}{\partial
\bkg\theta_k}}\big)^{T}. $ \item[(d)] For a given covariate $\bk
x$, we can estimate $\bkg\Sigma_{|\bk x}$ by
\begin{eqnarray*}
\widehat{\bkg\Sigma_{|\bk x}}&=
diag(\bkg\pi(\widehat{\bkg\theta}_n,\bk
x))-\bkg\pi(\widehat{\bkg\theta}_n,\bk
x)^{T}\bkg\pi(\widehat{\bkg\theta}_n,\bk x)~
~~~~~~~~~~~~~~~~~~~~~~~~~~~~~~~~~~~~ &\\
&+ 2\sum_{k=1}^K \Big(\frac{\partial \bkg\pi(\bkg\theta^{\ast},\bk
x)}{\partial
\bkg\theta_k^{\ast}}\Big|_{\smallbkg\theta^{\ast}=\widehat{\smallbkg\theta}_n}\Big)
\widehat{\bk V}_{n,k} \Big(\frac{\partial
\bkg\pi(\bkg\theta^{\ast},\bk x)}{\partial \bkg\theta_k^{\ast}}
\Big|_{\smallbkg\theta^{\ast}=\widehat{\smallbkg\theta}_n}\Big)^{T}
\frac{\#\{m\le n:\bkg\xi_m=\bk x\}}{n}.&
\end{eqnarray*}
\end{enumerate}
Notice that $\phi_k\bk I_k(\bkg\theta|\bkg\xi)$ does not depend on
$\phi_k$. When $\phi_k$ is unknown, we can estimate $\bk I_k$ in
the same way after replacing $\phi_k$ with its estimate
$\widehat{\phi}_k$.

\end{remark}

We now consider two examples, the logistic regression model and the
linear model.

\smallskip

{\it Example 3.1. Logistic Regression Model.}~~ We consider the case of
dichotomous (i.e., success or failure) responses. Let $Y_k=1$ if a
subject being given treatment $k$ is a success and $0$ otherwise,
$k=1,\ldots,K$. Let
$p_k=p_k(\bkg\theta_k,\bkg\xi)=\pr(Y_k=1|\bkg\xi)$ be the
probability of the success of a trial of treatment $k$ for a given
covariate $\bkg\xi$, $q_k=1-p_k$, $k=1,\ldots,K$. Assume that
\begin{eqnarray}\label{eqLogistM}
\logit(p_k)=\alpha_k+\bkg\theta_k\bkg\xi^T, \quad k=1,\ldots,K.
\end{eqnarray}
Without loss of generality, we assume that $\alpha_k=0$,
$k=1,2,\ldots,K$, or alternatively, we can redefine the covariate
vector to be $(1,\bkg\xi)$.  For each $k=1,\ldots,K$, let
$p_{j,k}=p_k(\bkg\theta_k,\bkg\xi_k)$. With the observations up to
stage $m$, the MLE $\widehat{\bkg\theta}_{m,k}$ of $\bkg\theta_k$
($k=1,\ldots,K$) is that for which $\widehat{\bkg\theta}_{m,k}$ maximizes
\begin{eqnarray}\label{MLE}
L_k=:\prod_{j=1}^m
p_{j,k}^{X_{j,k}Y_{j,k}}(1-p_{j,k})^{X_{j,k}(1-Y_{j,k})} \; \text{
over } \bkg\theta_k\in \bkg\Theta_k.
\end{eqnarray}
The logistic regression model is a special case of GLM (\ref{eqglm})
with $\phi_k=1$, $\mu_k=\log(p_k/q_k)$, $h_k(x)=x$,
$b_k(y_k,\phi_k)=0$, and
$a_k(\mu_k)=-\log(1-p_k)=\log(1+e^{\mu_k})$. Thus, given $\bkg\xi$,
the conditional information matrix is $\bk
I_k(\bkg\theta_k|\bkg\xi)=a_k^{\prime\prime}(\mu_k)\bkg\xi^{T}\bkg\xi
=p_kq_k\bkg\xi^{T}\bkg\xi$. For Theorem \ref{thglm}, we have the
following corollary.

\begin{corollary} \label{thlogistic}
 Suppose that Condition \ref{conditionA} is satisfied,
 $\ep\|\bkg\xi\|^2<\infty$, and the matrix $\ep[\bkg\xi^{T}\bkg\xi]$ is nonsingular.
We then have (\ref{eqConsiAssProb}), (\ref{eqLIL}), (\ref{eqCLT})
with $\bk V_k=\bk I_k^{-1}$ and
 $\bk I_k=\ep\{\pi_k(\bkg\theta,\bkg\xi)p_k q_k \bkg\xi^{T}\bkg\xi\}, \;\;
k=1,\ldots, K.
$
 Moreover, if $\pr(\bkg\xi=\bk x)>0$ for a given covariate $\bk x$,
 then (\ref{givenCovConsi}) and (\ref{givenCovNormal}) hold.
\end{corollary}

\smallskip
{\it Example 3.2. Normal Linear Regression Model.} The responses
are normally distributed, that is, $Y_k|_{\smallbkg\xi}\sim
N(\mu_k,\sigma_k^2)$ with link function
$\mu_k=\bkg\theta_k\bkg\xi^{T}$, then the linear model is a
special case of GLM (\ref{eqglm}) with $\phi_k=\sigma_k^2$,
$a_k(\mu_k)=\mu_k^2/2$ and $h_k(x)=x$. Thus, we have the following
corollary.

\begin{corollary}\label{thlm}
 Suppose that the conditions in Corollary \ref{thlogistic} are
 satisfied.
We then have (\ref{eqConsiAssProb}), (\ref{eqLIL}),
(\ref{eqCLT}) with
 $\bk V_k=\bk I_k^{-1}$ and
 $\bk I_k=\ep[\pi_k(\bkg\theta,\bkg\xi)\bkg\xi^{T}\bkg\xi]/\sigma_k^2$,
$k=1,\ldots,K$. Moreover, if $\pr(\bkg\xi=\bk x)>0$ for given $\bk
x$, then (\ref{givenCovConsi}) and (\ref{givenCovNormal}) hold.
\end{corollary}

\begin{remark}\label{rek3.0}
Bandyopadhyay and Biswas (2001) considered the normal linear
regression model in which $\theta_{11}=\mu_1$,
$\theta_{21}=\mu_2$, $\theta_{1j}=\theta_{2j}=\beta_{j-1}$,
$j=2,\ldots,d$, and the first component of $\bkg\xi$ is $1$. Their
proposed allocation probabilities are functions of estimates of
the unknown parameters that depend only on information of the
previous patients, but not on the covariates of the incoming
patient. Theorem 1 of Bandyopadhyay and Biswas (2001) gives the
consistency property of $N_{n,1} /n$ and $\pr(X_{n,1}=1)$.
However, their proof is not correct, since the assignments
$\delta_1,\cdots,\delta_i$ are functions of the previous responses
$Y_1,\cdots,Y_{i-1}$ and covariates.  In fact, given the
assignments $\delta_1,\cdots,\delta_i$, the responses
$Y_1,\cdots,Y_i$ are no longer independent normal variables, which
implies that their equation (4) is not valid.  Nevertheless, if we
let
$$\bkg\xi=(1,\widetilde{\bkg\xi}), ~\bk a=\ep
\widetilde{\bkg\xi},~ \bk
I_{\widetilde{\xi}}=\Var\{\widetilde{\bkg\xi}\},~
v_1=\Phi(\frac{\mu_1-\mu_2}{T}) \mbox{ and }v_2=1-v_1,$$ under our
theoretical framework, it can be proved that  Theorem 1 of
Bandyopadhyay and Biswas (2001) is correct.  Further, it is not
difficult to show that
$$ \sqrt{n}( \widehat{\mu}_{n1}-\mu_1,
\widehat{\mu}_{n2}-\mu_2)\overset{D}\to N\left((0,0),
\sigma^2\begin{pmatrix} 1/v_1+\bk a\bk
I_{\widetilde{\xi}}^{-1}\bk a^T & \bk a\bk I_{\widetilde{\xi}}^{-1}\bk a^T \\
\;\; \bk a\bk I_{\widetilde{\xi}}^{-1}\bk a^T  & 1/v_2+\bk a\bk
I_{\widetilde{\xi}}^{-1}\bk a^T
\end{pmatrix}\right),
$$
$$ \sqrt{n}(\widehat{\beta}_{n}-\beta)\overset{D}\to N(\bk 0, \sigma^2\bk I_{\widetilde{\xi}}^{-1}), $$
and
$$ \sqrt{n}\left( N_{n1}/n-v_1\right)\overset{D}\to N\left(0,v_1v_2+
\frac{2\sigma^2}{v_1v_2}\left(\Phi^{\prime}(\frac{\mu_1-\mu_2}{T})\right)^2\right),
$$
where $\sigma^2$ is the variance of the errors in the linear model.
\end{remark}

\begin{remark}\label{rek3.1}
Corollary \ref{thlm} can be generalized to general responses.
Suppose that the response
$Y_k$ of a subject to treatment $k$, $k=1,\ldots,K$,   and its
covariate $\bkg\xi$ satisfies the linear regression model
\begin{eqnarray*}
\ep[Y_k|\bkg\xi]=p_k(\bkg\theta_k,\bkg\xi)=\bkg\theta_k\bkg\xi^{T},
\quad k=1,\ldots, K.
\end{eqnarray*}
For the observations up to stage $m$,  let
$\widehat{\bkg\theta}_{m,k}$ minimize the error sum of squares
\begin{eqnarray*}
S_k(\bkg\theta_k) = \sum_{j=1}^m
X_{j,k}(Y_{j,k}-\bkg\theta_k\bkg\xi_j^{T})^2 \text{ over }
\bkg\theta_k \in \bkg\Theta_k,
\end{eqnarray*}
$k=1,\ldots, K$. Here, $\widehat{\bkg\theta}_{m,k}$ is the
least-squares estimator (LSE) of $\bkg\theta_k$. Then Corollary
\ref{thlm} remains true with
 $\bk V_k=\bk I_{\xi k}^{-1}\bk I_{Y k} \bk
I_{\xi k }^{-1}$, $
 \bk I_{\xi k
}=\ep[\pi_k(\bkg\theta,\bkg\xi)\bkg\xi^{T}\bkg\xi]$ and $\bk I_{Y k
} =\ep\{\pi_k(\bkg\theta,\bkg\xi)(Y_k-\bkg\theta_k\bkg\xi^{T})^2
\bkg\xi^{T}\bkg\xi\}$ under the condition
$\ep\|Y_k\bkg\xi\|^2<\infty$,  $k=1,\ldots,K$. This result follows
from Theorem \ref{th2}, as condition (\ref{rateconsi}) is satisfied
with $\bk h_k=(Y_k-\bkg\theta_k\bkg\xi^{T})\bkg\xi \bk I_{\xi
k}^{-1}$.
\end{remark}

\smallskip
\setcounter{section}{5} \setcounter{equation}{0}
\setcounter{theorem}{0} \setcounter{remark}{0}
{\bf 4. Discussion.}

This paper makes two major contributions.  First,
a comprehensive framework of  CARA
designs is proposed to serve as a paradigm for treatment
allocation procedures in clinical trials when covariates are
available. It is a very general framework that allows a wide
spectrum of applications to very general statistical models,
including generalized linear models as special cases. Second,
asymptotic properties are obtained to provide a statistical basis
for inferences after using a CARA design.

When covariate information is not being used in the treatment
allocation scheme, an optimal allocation proportion is usually
determined with the assistance of some optimality criteria.
Jennison and Turnbull (2000) described a general procedure to
search for an optimal allocation. For CARA
designs, the means to
define and obtain an optimal allocation scheme is still unclear.
For CARA design, we can find optimal allocation for each fixed value of the covariate.
Theorem 2.2 provides theoretical support for targeting optimal
allocation by using CARA design for each fixed covariate.

For response-adaptive designs without covariates, Hu and
Rosenberger (2003) studied the relationship among the power, the
target allocation, and the variability of the designs. It is
important to study the behavior of the power function when a
CARA design is used in clinical
trials. It is not
difficult to derive the power function for binary responses with discrete
covariates. For the general
covariate $\bkg\xi$, the formulation becomes very different, and it
is an interesting topic for future research.

\renewcommand{\theequation}{A.\arabic{equation}}
\setcounter{equation}{0}
\newcommand{\proofend}{$\quad \Box$}
\begin{center}
{APPENDIX: PROOFS}
\end{center}

The proofs of the theorems are organized as follows. First,
we prove the theorems for the general CARA
design in Section 2. We then derive the results in
Section 3 by the application of the theorems in Section 2.

\smallskip

\noindent
{\sc Proof of Theorem \ref{th2}.} First, notice that for each
$k=1,\ldots, K$,
$ X_{m+1,k}=X_{m+1,k}-\ep[X_{m+1,k}|\Cal
F_m]+g_k(\widehat{\bkg\theta}_m) $ and then
\begin{eqnarray}\label{eqproof1.1}
N_{n,k}=\ep[X_{1,k}|\Cal F_0]+\sum_{m=1}^n
(X_{m,k}-\ep[X_{m,k}|\Cal
F_{m-1}])+\sum_{m=1}^{n-1}g_k(\widehat{\bkg\theta}_m).
\end{eqnarray}
The second term is a martingale. We then show that the third
term can be approximated by another martingale.  Write $\Delta
M_{m,k}=X_{m,k}-\ep[X_{m,k}|\Cal F_{m-1}]$, $\Delta \bk
T_{m,k}=X_{m,k}\bk h_k(Y_{m,k},\bkg\xi_m)$, $k=1,\ldots, K$. Let $
\bk M_n=\sum_{m=1}^n \Delta \bk M_m$ and $ \bk T_n=\sum_{m=1}^n
\Delta \bk T_m$, where $\Delta \bk M_m=(\Delta M_{m,1},\cdots, $
$\Delta M_{m,K})$ and $\Delta \bk T_m=(\Delta \bk T_{m,1},\cdots,
\Delta \bk T_{m,K})$. Here, the symbol $\Delta$ denotes the
differencing operand of a sequence $\{z_n\}$, i.e., $\Delta
z_n=z_n-z_{n-1}$.  Then $\{(\bk M_n, \bk T_n)\}$ is a
multi-dimensional martingale sequence that satisfies
\begin{eqnarray}\label{eqproofmomentM}
\quad \quad \quad |\Delta M_{n,k}|\le 1,\quad
 \|\Delta \bk T_{n,k}\|\le \| \bk
h_k(Y_{n,k},\bkg\xi_n)\|
\end{eqnarray}
and $\ep\| \bk h_k(Y_{n,k},\bkg\xi_n)\|^2<\infty$, $k=1,\ldots, K$.
It follows that
\begin{eqnarray}\label{eqproofML2}
\|\bk M_n\|=O(\sqrt{n})\;\; \mbox{and} \;\; \|\bk T_n\|=O(\sqrt{n})
\;\;\text{ in }\; L_2.
\end{eqnarray}
Also, according to the law of the iterated logarithm for
martingales, we have
\begin{eqnarray} \label{LILforMartingale}
\bk M_n =O\big(\sqrt{n\log\log n}\big)\;\; a.s. \text{ and }
\bk T_n =O\big(\sqrt{n\log\log n}\big)\;\; a.s..
\end{eqnarray}
 From
(\ref{LILforMartingale}) and (\ref{rateconsi}), it follows that
\begin{eqnarray}\label{forthetaLIL}
\widehat{\bkg\theta}_n-\bkg\theta=O\Big(\sqrt{\frac{\log\log
n}{n}}\Big)\;\; a.s.
\end{eqnarray}
From (\ref{eqproof1.1}), (\ref{forthetaLIL}),
(\ref{eqproofML2}), and Condition \ref{conditionA}, it follows that
\begin{eqnarray*}
N_{n,k}-n v_k&=&M_{n,k}+\sum_{m=1}^{n-1}\sum_{j=1}^K
(\widehat{\bkg\theta}_{m,j}-\bkg\theta_j)\big(\frac{\partial
g_k}{\partial \bkg\theta_j }\big)^{T}
+\sum_{m=1}^{n-1}o(\|\widehat{\bkg\theta}_m-\bkg\theta\|^{1+\delta})
\nonumber\\
&=&M_{n,k}+\sum_{m=1}^n \sum_{j=1}^K\frac{\bk
T_{m,j}\big(1+o(1)\big)}{m} \big(\frac{\partial g_k}{\partial
\bkg\theta_j }\big)^{T}+o(n^{1/2})\;\; a.s..\\
&=&M_{n,k}+\sum_{m=1}^n \sum_{j=1}^K\frac{\bk T_{m,j}}{m}
\big(\frac{\partial g_k}{\partial \bkg\theta_j
}\big)^{T}+o(n^{1/2})\;\;\text{ in probability},
\end{eqnarray*}
that is,
\begin{eqnarray}\label{eqproofth2.5}
 ~~~~~\bk N_n-n \bk v &=& \bk M_n+ \sum_{m=1}^n \sum_{j=1}^K\frac{\bk
T_{m,j}\big(1+o(1)\big)}{m} \big(\frac{\partial \bk g}{\partial
\bkg\theta_j }\big)^{T}+o(n^{1/2})\;\;a.s.
\\
\label{eqproofth2.5ad} &=& \bk G_n+o(n^{1/2})\;\;\text{ in
probability},
\end{eqnarray}
where
$$ \bk G_n=\bk M_n+ \sum_{m=1}^n \sum_{j=1}^K\frac{\bk
T_{m,j}}{m} \big(\frac{\partial \bk g}{\partial \bkg\theta_j
}\big)^{T}=\bk M_n+ \sum_{m=1}^n \sum_{j=1}^K\Delta\bk T_{m,j}
\big(\frac{\partial \bk g}{\partial \bkg\theta_j
}\big)^{T}\sum_{i=m}^n\frac{1}{i}.
$$
The combination of (\ref{LILforMartingale}) and (\ref{eqproofth2.5}) yields
$$ \bk N_n-n \bk v= O\big(\sqrt{n\log\log n}\big)
+\sum_{m=1}^n \sum_{j=1}^K\frac{O(\sqrt{m\log\log m})}{m}=
O\big(\sqrt{n\log\log n}\big) \;\; a.s. $$ (\ref{eqConsiAssProb})
is obvious by noting (\ref{forthetaLIL}) and the continuity of
$\bkg\pi(\cdot,\bkg\xi)$. The proof of consistency is thus obtained.
Next, we consider the asymptotic normality. Notice that $\bk M_n$, $\bk
T_n$, and $\bk G_n$ are all the sum of martingale differences. It
is easy to verify that  the Lindberg condition is satisfied by
(\ref{eqproofmomentM}). To complete the proof it suffices to
derive the variances. First, the conditional variance-covariance
matrices of the martingale difference $\{\Delta\bk M_n,\Delta\bk
T_n\}$ satisfy
$$ \ep[(\Delta\bk M_n)^{T}\Delta\bk M_n|\Cal F_{n-1}]
=diag\big(\bk g(\widehat{\bkg\theta}_{n-1})\big)-\big(\bk
g(\widehat{\bkg\theta}_{n-1})\big)^{T}\bk
g(\widehat{\bkg\theta}_{n-1}) \to \bkg\Sigma_1 \; \; \text{ in }
L_1;
$$
$$ \ep[(\Delta\bk T_{n,k})^{T}\Delta\bk T_{n,k}|\Cal F_{n-1}]
=\ep\big[\pi_k(\bk x,\bkg\xi)\big(\bk h_k(Y_k,\bkg\xi)\big)^{T}\bk
h_k(Y_k,\bkg\xi)\big]\big|_{\bk x=\widehat{\smallbkg\theta}_{n-1}}
\to \bk V_k \; \; \text{ in } L_1;
$$
$$ \ep[(\Delta M_{n,i})^{T}\Delta\bk T_{n,j}|\Cal F_{n-1}]=\bk
0 \; \text{ for all } i, j\; \;\text{ and } \; \ep[(\Delta\bk
T_{n,i})^{T}\Delta\bk T_{n,j}|\Cal F_{n-1}]=\bk 0\;
$$
for all $i\ne j$.  It follows that $\Var\{\bk T_n\}/n \to \bk V$ and
\begin{eqnarray*}
\Var\{\bk
G_n\}&=&\sum_{m=1}^n[\bkg\Sigma_1+o(1)]+\sum_{m=1}^n\sum_{l=1}^n
\sum_{j=1}^K \frac{l\wedge m}{m l} \frac{\partial \bk g}{\partial
\bkg\theta_j }[\bk V_j+o(1)] \big(\frac{\partial \bk g}{\partial
\bkg\theta_j }\big)^{T}\\
&=&n(\bkg\Sigma_1+2\bkg\Sigma_2)+o(n)=n\bkg\Sigma+o(n).
\end{eqnarray*}
By the central limit theorem for martingales (Hall and
Heyde, 1980), it follows that
$
\sqrt{n}(\widehat{\bkg\theta}_n-\bkg\theta)=n^{-1/2}\bk
T_n+o(1)\overset{\Cal D}\to N(\bk 0,\bk V)
$
and
\begin{eqnarray}\label{eqforN}
\sqrt{n}(\bk N_n/n-\bk v)=n^{-1/2}\bk G_n+o(1)\overset{\Cal D}\to
N(\bk 0,\bkg\Sigma).
\end{eqnarray}
The proof is now complete. \proofend

\smallskip
\noindent
{\sc Proof of Theorem \ref{givenCov}.} First, according to the law of
large numbers, we have
\begin{eqnarray}\label{eqproofth2.7} \frac 1n \sum_{m=1}^n I\{\bkg \xi_m=\bk x\}
\to \pr(\bkg\xi=\bk x)\;\; a.s. \end{eqnarray}
 and
\begin{eqnarray*}
\frac 1n \sum_{m=1}^n X_{m,k}I\{\bkg \xi_m=\bk x\}
&=& \frac 1n \sum_{m=1}^n \Big(X_{m,k}I\{\bkg \xi_m=\bk x\}-
\ep[X_{m,k}I\{\bkg \xi_m=\bk x\}|\Cal F_{m-1}]\Big)\\
&&+\frac 1n \sum_{m=1}^n \pi_k(\widehat{\bkg\theta}_{m-1}, \bk
x)\pr(\bkg\xi_m=\bk x) \to \pi_k(\bkg\theta, \bk x)\pr(\bkg\xi=\bk
x) \;\; a.s.,
\end{eqnarray*}
and thus (\ref{givenCovConsi}) is proved. We then consider the asymptotic
 normality. The proof is similar to that of
 (\ref{eqforN}). The difference lies in the approximation of the
 process by a new $2K$ dimensional martingale and the calculation of its
  variance-covariance matrix. Define
 $\zeta_{n,k}(\bk x):=\sum_{m=1}^n (X_{m,k}-\pi_k(\bkg\theta,\bk x))
I\{\bkg\xi_m=\bk x\}.$ Then,
$$
 \sqrt{N_n(\bk x)}\Big(\frac{\bk N_{n|\bk x}}{N_n(\bk x)}-
\bkg\pi(\bkg\theta,\bk x)\Big)=\sqrt{\frac{n}{N_n(\bk x)}}
\frac{\zeta_{n,k}(\bk x)}{\sqrt{n}}, \; \; k=1,\ldots, K.
$$
Notice (\ref{eqproofth2.7}). It is sufficient to prove
\begin{eqnarray}\label{eqproofgiven1.1}
n^{-1/2}\big(\zeta_{n,1}(\bk x),\ldots,\zeta_{n,K}(\bk
x)\big)\overset{\Cal D}\to N\big(\bk 0,\bkg\Sigma_{|\bk x}
\pr(\bkg\xi=\bk x)\big).
\end{eqnarray}
 With the same argument as is used to derive
(\ref{eqproofth2.5}), we can obtain
\begin{eqnarray*}
  \zeta_{n,k}(\bk x) &=&\sum_{m=1}^n \big(\Delta \zeta_{n,k}(\bk x)-
\ep[\Delta\zeta_{n,k}(\bk x)|\Cal
F_{n-1}]\big) \\
&& +\sum_{m=1}^n\big(\pi_k(\widehat{\bkg\theta}_{m-1}, \bk x)-\pi_k(\bkg\theta,
\bk x)\big)\pr(\bkg\xi=\bk x)\\
&=&\sum_{m=1}^n \big(\Delta \zeta_{m,k}(\bk x)-
\ep[\Delta\zeta_{m,k}(\bk x)|\Cal F_{m-1}]\big)\\
&&+\sum_{j=1}^K\sum_{m=1}^n\frac {\bk
T_{m,j}}{m}\Big(\frac{\partial \pi_k(\bkg\theta,\bk
x)}{\partial\bkg\theta_j}\Big)^{T}\pr(\bkg\xi=\bk x)+o(n^{1/2})
\;\; \text{in probability}.\\
&=&:G_{n,k}(\bk x)+o(n^{1/2}).
\end{eqnarray*}
Similar to the proof of (\ref{eqforN}), to complete the proof it suffices
to get the variance of $\bk G_n(\bk x)=(G_{n,1}(\bk x),
\ldots, G_{n,K}(\bk x))$. Let $\Delta \overline {M}_{n,k}(\bk
x)=\Delta \zeta_{n,k}(\bk x)-\ep[\Delta\zeta_{n,k}(\bk x)|\Cal
F_{n-1}]$. The variance-covariance matrix of the martingale
difference $\{(\Delta\overline{\bk M}_n(\bk x), \Delta\bk T_n)\}$ then
satisfies
$\ep[ (\Delta \overline{M}_{n,k}(\bk x))^2|\Cal F_{n-1}]\to
\pi_k(\bkg\theta,\bk x)\big(1-\pi_k(\bkg\theta,\bk
x)\big)\pr(\bkg\xi=\bk x)$,

\noindent
$\ep[ \Delta
\overline{M}_{n,k}(\bk x) \Delta \overline{M}_{n,j}(\bk x)|\Cal
F_{n-1}] \to -\pi_k(\bkg\theta,\bk x)\pi_j(\bkg\theta,\bk
x)\pr(\bkg\xi=\bk x)\;\forall k\ne j,
$ and

\noindent
$\ep[ \Delta \overline{M}_{n,k}(\bk x) \Delta\bk T_{n,j}|\Cal F_{n-1}]=\bk
0\;  \forall i,j
$
in $L_1$. It follows that
\begin{eqnarray*}
\Var\{\bk G_n(\bk x)\} &=&n\left[diag(\bkg\pi(\bkg\theta,\bk
x)-\bkg\pi(\bkg\theta,\bk x)^T\bkg\pi(\bkg\theta,\bk
x)\pr(\bkg\xi=\bk x)+o(1)\right]
\\
&&+ \sum_{m=1}^n\sum_{l=1}^n\sum_{j=1}^K\frac{l\wedge m}{l m}
\frac{\partial \bkg\pi(\bkg\theta,\bk x)}{\partial
\bkg\theta_j}[\bk V_j+o(1)] \big(\frac{\partial
\bkg\pi(\bkg\theta,\bk x)}{\partial \bkg\theta_j}\big)^{T}
\pr^2(\bkg\xi=\bk x)
\\
&=& n\bkg\Sigma_{|\bk x}\pr(\bkg\xi=\bk x)+o(n).
\end{eqnarray*}
 (\ref{eqproofgiven1.1}) is then proved.  \proofend

\smallskip
\noindent
{\sc Proof of Corollary \ref{thglm}.}  By Theorem \ref{th2}, it
suffices to verify the condition (\ref{rateconsi}). Notice that
$\widehat{\bkg\theta}_{m,k}$ is a solution to $\partial \log
L_k/\partial \bkg\theta_k=0$. The application of Taylor's theorem yields
\begin{eqnarray}\label{eqproofthglm1.1}
\frac{\partial \log L_k}{\partial
\bkg\theta_k}\Big|_{\smallbkg\theta_k}
+(\widehat{\bkg\theta}_{m,k}-\bkg\theta_k)\Big\{\frac{\partial^2
\log L_k}{\partial \bkg\theta_k^2}\Big|_{\smallbkg\theta_k}
+\int_0^1\Big[\left.\frac{\partial^2 \log L_k}{\partial
\bkg\theta_k^2}
\right|_{\smallbkg\theta_k}^{\smallbkg\theta_k+t(\widehat{\smallbkg\theta}_{m,k}-
\smallbkg\theta_k)}\Big]dt\Big\}
\nonumber\\
=\frac{\partial \log L_k}{\partial
\bkg\theta_k}\Big|_{\widehat{\smallbkg\theta}_{m,k}}=\bk 0,
\end{eqnarray}
where $f(x)\Big|_{a}^b=f(b)-f(a)$. Notice that
\begin{eqnarray}\label{eqproofthglm1.2}
\frac{\partial \log L_k}{\partial \bkg\theta_k}=\sum_{j=1}^m
X_{j,k}\frac{\partial\log
f_k(Y_{j,k}|\bkg\xi_j,\bkg\theta_k)}{\partial \bkg\theta_k}
\end{eqnarray}
and
$$\frac{\partial^2
\log L_k}{\partial \bkg\theta_k^2}=\sum_{j=1}^m
X_{j,k}\frac{\partial^2\log
f_k(Y_{j,k}|\bkg\xi_j,\bkg\theta_k)}{\partial \bkg\theta_k^2}.
$$
We assume that the following regular condition
\begin{eqnarray}\label{eqregular}
H(\delta)=:\ep\Bigl[\sup_{\|\bk z\|\le
\delta}\Bigl\|\frac{\partial^2\log
f_k(Y_k|\bkg\xi,\bkg\theta_k)}{\partial \bkg\theta_k^2}
\biggr|_{\smallbkg\theta_k}^{\smallbkg\theta_k+\bk
z}\Bigr\|\Bigl]\to 0 \text{ as } \delta\to 0.
\end{eqnarray}
This regularity condition is implied by the simple condition that
$a_k^{\prime\prime}$, $h_k^{\prime\prime}$ are continuous and
$\bkg\xi$ is bounded.  Under (\ref{eqregular}), one can show that
$$\sup_{\|\bk z\|\le \delta}\Bigl\|\frac{1}{m}\frac{\partial^2 \log L_k}{\partial
\bkg\theta_k^2} \biggr|_{\smallbkg\theta_k}^{\smallbkg\theta_k+\bk
z}\Bigr\|\le H(\delta)+o(1)\;\; a.s..
$$
However,
$$\sum_{j=1}^m\left\{
X_{j,k}\frac{\partial^2\log
f_k(Y_{j,k}|\bkg\xi_j,\bkg\theta_k)}{\partial
\bkg\theta_k^2}-\ep\Big[X_{j,k}\frac{\partial^2\log
f_k(Y_{j,k}|\bkg\xi_j,\bkg\theta_k)}{\partial
\bkg\theta_k^2}\Big|\Cal F_{j-1}\Big]\right\}
$$
is a martingale. According to the law of large numbers,
\begin{eqnarray}\label{eqproofthglm1.3}\frac{\partial^2 \log L_k}{\partial
\bkg\theta_k^2}&=&\sum_{j=1}^m\ep\Big[X_{j,k}\frac{\partial^2\log
f_k(Y_{j,k}|\bkg\xi_j,\bkg\theta_k)}{\partial
\bkg\theta_k^2}\Big|\Cal
F_{j-1}\Big]+o(m)\nonumber\\
&=&-\sum_{j=1}^m\big\{\ep\big[\pi_k(\bk z,\bkg\xi)\bk
I_k(\bkg\theta_k|\bkg\xi)\big]\big\}_{\bk
z=\widehat{\smallbkg\theta}_{j-1}}+o(m) = -m\bk I_k+o(m)\;\; a.s.
\end{eqnarray}
The substitution of (\ref{eqproofthglm1.2}) and  (\ref{eqproofthglm1.3})
into (\ref{eqproofthglm1.1}) yields
\begin{eqnarray*}
m(\widehat{\bkg\theta}_{m,k}-\bkg\theta_k)\Big\{\bk
I_k+o(1)+O\big(H(\|\widehat{\bkg\theta}_{m,k}-\bkg\theta_k\|)\big)\big\}
=\sum_{j=1}^mX_{j,k}\frac{\partial\log
f_k(Y_{j,k}|\bkg\xi_j,\bkg\theta_k)}{\partial \bkg\theta_k}.
\end{eqnarray*}
Thus,
\begin{eqnarray}\label{appthetainglm}
~~~~~~~~\widehat{\bkg\theta}_{m,k}-\bkg\theta_k
=\frac{1}{m}\sum_{j=1}^mX_{j,k}\frac{\partial\log
f_k(Y_{j,k}|\bkg\xi_j,\bkg\theta_k)}{\partial \bkg\theta_k}\bk
I_k^{-1}\big(1+o(1)\big)\;\; a.s.
\end{eqnarray}
Notice that
$$\ep[\frac{\partial\log
f_k(Y_{j,k}|\bkg\xi_j,\bkg\theta_k)}{\partial
\bkg\theta_k}|\bkg\xi_j]=0 \;\;\text{ and }\;\;
\Var\big\{\frac{\partial\log
f_k(Y_{j,k}|\bkg\xi_j,\bkg\theta_k)}{\partial
\bkg\theta_k}|\bkg\xi_j\big\}=\bk I_k(\bkg\theta_k|\bkg\xi_j).$$
Hence, Condition (\ref{rateconsi}) is valid. By Theorem \ref{th2},
the proof is complete.
\proofend

\begin{center} ACKNOWLEDGEMENTS \end{center}
\vspace{-0.1in} Special thanks go to the anonymous referees, the associate editor and
the editor for their constructive comments,
which led to a much improved
version of the paper.

\smallskip

\baselineskip 14pt
\begin{center}
REFERENCES
\end{center}

\vspace{-0.1in}
\footnotesize

\begin{enumerate}
\item[{[1]}] {\sc Atkinson,} A. C. (1982). Optimal biased coin designs
for sequential clinical trials with prognostic factors. {\em
Biometrika} {\bf 69} 61-67.
\item[{[2]}] {\sc Atkinson,} A. C. (1999).
 Optimal biased-coin designs for
sequential treatment allocation with covariate information.
{\em Statist. Med.}  {\bf 18} 1741-1752.

\item[{[3]}] {\sc Atkinson,} A. C. (2002). The comparison of designs for
sequential clinical trials with covariate information. {\em
Journal of the Royal Statistical Society - Series A}  {\bf 165}
349-373.

\item[{[4]}] {\sc Atkinson,} A. C. (2004). Adaptive biased-coin designs
for clinical trials with several treatments. {\em Discussiones
Mathematicae Probability and Statistics} {\bf 24} 85-108.

\item[{[5]}] {\sc Atkinson,} A. C. and {\sc Biswas,} A. (2005a).
Adaptive biased-coin designs for skewing the allocation proportion in clinical trials with normal
responses. {\em Statistics in Medicine} {\bf 24} 2477-2492.

\item[{[6]}] {\sc Atkinson,} A. C. and {\sc Biswas,} A. (2005b).
Bayesian adaptive biased-coin designs for clinical trials with normal
responses. {\em Biometrics} {\bf 61} 118-125.

\item[{[7]}] {\sc Baldi Antognini}, A. and {\sc Giovagnoli}, A. (2004).
On the large sample optimality of sequential designs for comparing
two-treatments. Manuscript, University of Bologna.

\item[{[8]}] {\sc Bandyopadhyay}, U. and {\sc Biswas,} A. (1999).
Allocation by randomized play-the-winner rule in the presence of prognostic factors.
{\em Sankhya} B {\bf 61} 397-412.

\item[{[9]}] {\sc Bandyopadhyay}, U. and {\sc Biswas,} A. (2001).
Adaptive designs for normal responses with prognostic factors.
{\em Biometrika} {\bf 88} 409-419.

\item[{[10]}] {\sc Eisele,} J. and {\sc Woodroofe,} M. (1995).
Central limit theorems for doubly adaptive biased coin designs.
{\em Ann. Statist.} {\bf 23} 234-254.

\item[{[11]}]  {\sc Efron,} B. (1971). Forcing a sequential
experiment to be balanced.
 {\em Biometrika} {\bf 62}  347-352.

\item[{[12]}] {\sc Hall, P.} and {\sc Heyde,} C. C. (1980). {\it
Martingale Limit Theory and its Applications}. Academic Press,
London.

\item[{[13]}] {\sc Hu,} F. and  {\sc Rosenberger,} W. F.  (2003).
Optimality, variability, power   evaluating response-adaptive
randomization procedures for treatment comparisons.  {\em J. Amer.
Statist. Assoc.} {\bf 98} 671-678.

\item[{[14]}] {\sc Jennison,} C. and {\sc Turnbull,} B. W. (2000).
{\em Group Sequential Methods with Applications to Clinical
Trials}. Chapman and Hall/CRC, Boca Raton, FL.

\item[{[15]}] {\sc McCullagh,} P. and {\sc Nelder,} J. A. (1989).
{\em Generalized Linear Models}, Second Edition, Chapman and Hall,
London.

\item[{[16]}] {\sc Pocock,} S. J. and {\sc Simon,} R. (1975).
Sequential treatment assignment with balancing for prognostic
factors in the controlled clinical trial.  {\em Biometrics} {\bf
31} 103-115.

\item[{[17]}] {\sc Robbins,} H. (1952). Some aspects of the
sequential design of experiments. {\em Bull. Amer. Math. Soc.}
{\bf 58} 527-535.

\item[{[18]}] {\sc Rosenberger,} W. F.,  {\sc Flournoy,} N. and
{\sc Durham}, S. D. (1997). Asymptotic normality of maximum
likelihood estimators from multiparameter response-driven design.
{\em J. Statist. Plann. Inf.} {\bf 60} 69-76.

\item[{[19]}] {\sc Rosenberger,} W. F. and {\sc Hu,} F. (2004).
Maximizing power and minimizing treatment failures. {\em Clinical
Trials}  {\bf 1} 141-147.

\item[{[20]}] {\sc Rosenberger,} W. F. and {\sc Hu,} M. X. (2002).
On the use of generalized linear models following a sequential
design. {\em Statist. Probab. Letters} {\bf 56} 155-161.

\item[{[21]}] {\sc Rosenberger,} W. F., {\sc Vidyashankar,} A. N.
and {\sc Agarwal,} D. K. (2001). Covariate-adjusted
response-adaptive designs for binary response. {\em J. Biopharm.
Statist.} {\bf 11} 227-236.

\item[{[22]}] {\sc Thompson,} W. R. (1933). On the likelihood that
one unknown probability exceeds another in view of the evidence of
the two samples. {\em Biometrika} {\bf 25} 275-294.

\item[{[23]}] {\sc Wei, L. J.} and {\sc Durham,} S. (1978). The
randomized pay-the-winner rule in medical trials. {\em J. Amer.
Statist. Assoc.} {\bf 73} 840-843.

\item[{[24]}] {\sc Zelen,} M.  (1969).  Play-the-winner rule and
the controlled clinical trial.  {\em J. Amer. Statist. Assoc.}
{\bf 64} 131-146.

\item[{[25]}] {\sc Zelen,} M.  (1974). The randomization and
stratification of patients to clinical trials. {\em Journal of
Chronic Diseases} {\bf 28} 365-375. \item[{[29]}] {\sc Zelen,} M.
and {\sc Wei,} L. J. (1995).  Foreword.  In {\em Adaptive Designs}
(N. Flournoy and W. F. Rosenberger, eds).  IMS, Hayward, CA.
\end{enumerate}

\newpage

\renewcommand{\theequation}{B.\arabic{equation}}
\setcounter{equation}{0}
\begin{center}
{APPENDIX B: Additional PROOFS}
\end{center}

{\bf Proof of the existence and consistency of the solution of
(\ref{eqproofthglm1.1}):} It suffices to show that, for any
$\delta>0$ small enough, with probability one for $m$ large enough we have
 \begin{equation}\label{eqB1} \log
L_k(\bkg\theta_k^{\ast})<\log L_k(\bkg\theta_k), \text{ if }
\|\bkg\theta_k^{\ast}-\bkg\theta_k\|=\delta.
\end{equation}
The application of Taylor's theorem yields
\begin{align*}
& \frac{1}{m}\log L_k(\bkg\theta_k^{\ast})-\frac{1}{m}\log
L_k(\bkg\theta_k)\nonumber\\
=&(\bkg\theta_k^{\ast}-\bkg\theta_k)\frac{1}{m}\frac{\partial \log
L_k}{\partial \bkg\theta_k}\Big|_{\smallbkg\theta_k}
+(\bkg\theta_k^{\ast}-\bkg\theta_k)\frac{1}{m}\frac{\partial^2 \log
L_k}{\partial
\bkg\theta_k^2}\Big|_{\smallbkg\theta_k}(\bkg\theta_k^{\ast}-\bkg\theta_k)^T
\nonumber\\
&
+(\bkg\theta_k^{\ast}-\bkg\theta_k)\Big\{\frac{1}{m}\int_0^1\Big[\left.\frac{\partial^2
\log L_k}{\partial \bkg\theta_k^2}
\right|_{\smallbkg\theta_k}^{\smallbkg\theta_k+t(\smallbkg\theta_k^{\ast}-
\smallbkg\theta_k)}\Big]dt\Big\}(\bkg\theta_k^{\ast}-\bkg\theta_k)^T.
\end{align*}
So  with probability one for $m$ large enough,
\begin{align*}
& \frac{1}{m}\log L_k(\bkg\theta_k^{\ast})-\frac{1}{m}\log
L_k(\bkg\theta_k)\nonumber\\
 \le &
-(\bkg\theta_k^{\ast}-\bkg\theta_k)\Big\{\frac{1}{m}\sum_{j=1}^m\big\{\ep\big[\pi_k(\bk
z,\bkg\xi)\bk I_k(\bkg\theta_k|\bkg\xi)\big]\big\}_{\bk
z=\widehat{\smallbkg\theta}_{j-1}}\Big\}(\bkg\theta_k^{\ast}-\bkg\theta_k)^T\\
&+\|\bkg\theta_k^{\ast}-\bkg\theta_k\|^2
H(\|\bkg\theta_k^{\ast}-\bkg\theta_k\|)+o(1)\nonumber\\
\le & -\|\bkg\theta_k^{\ast}-\bkg\theta_k\|^2 \min_{\|\bk y\|=1, \bk z\in
\smallbkg\Theta_k}\Big\{\bk y\big\{\ep\big[\pi_k(\bk z,\bkg\xi)\bk
I_k(\bkg\theta_k|\bkg\xi)\big]\big\}\bk y^T\Big\}+\|\bkg\theta_k^{\ast}-\bkg\theta_k\|^2
H(\|\bkg\theta_k^{\ast}-\bkg\theta_k\|)+o(1)\nonumber\\
\le& -c_0\delta^2+\delta^2H(\delta)+o(1)<0\; \text{uniformly in }
\bkg\theta_k^{\ast} \text{ with }
\|\bkg\theta_k^{\ast}-\bkg\theta_k\|=\delta
\end{align*}
 when $\delta$ is small enough. (\ref{eqB1}) is proved. \proofend

 \bigskip

 {\bf Proof of Corollary \ref{thlogistic}:} Notice
 $$ \frac{\partial^2\log
f_k(Y_k|\bkg\xi,\bkg\theta_k)}{\partial \bkg\theta_k^2}
=-p_kq_k\bkg\xi^T\bkg\xi $$
is bounded by $\|\bkg\xi\|^2$ and is continuous in $\bkg\theta_k$.
It follows that the regularity condition (\ref{eqregular}) is satisfied due to the domained convergence theorem.
On the other hand, it is obviously that
$$ \frac{\partial^2\log L_k}{\partial \bkg\theta^2}=-\sum_{j=1}^m X_{j,k}
p_k(\bkg\theta_k,\bkg\xi_j)q_k(\bkg\theta_k,\bkg\xi_j)\bkg\xi^T\bkg\xi$$
is a negatively definite matrix, and so $\log L_k(\bkg\theta_k^{\ast})$ is a strictly concave function of $\bkg\theta_k^{\ast}$.
It follows that the MLE is unique. Corollary \ref{thlogistic} now follows from Corollary \ref{thglm}.\proofend

\bigskip
{\bf Proof of Remark \ref{rek3.1}:}
It is obviously that $S_k(\bkg\theta_k)$ is strictly convex function of $\bkg\theta_k$.
It follows that the LSE $\widehat{\bkg\theta}_{m,k}$ exists and is unique.
On the hand, it is easily seen that
$\widehat{\bkg\theta}_{m,k}$ is the solution of the normal
equation as
\begin{equation}\label{eqB3}
\big(\widehat{\bkg\theta}_{m,k}-\bkg\theta_k\big)\frac{1}{m}
[\sum_{j=1}^m X_{j,k}\bkg\xi_j^{T}\bkg\xi_j]
=\frac{1}{m}\sum_{j=1}^m
X_{j,k}(Y_{j,k}-\bkg\theta_k\bkg\xi_j^{T})\bkg\xi_j.
\end{equation}
Also,
$\{X_{j,k}\bkg\xi_j^{T}\bkg\xi_j-\ep[X_{j,k}\bkg\xi_j^{T}\bkg\xi_j|
 \Cal F_{j-1}] \}$ and $\{X_{j,k}(Y_{j,k}-\bkg\theta_k\bkg\xi_j^{T})\bkg\xi_j\}$
  are both  sequences of martingale differences. It
follows from the law of large numbers for martingales that,
\begin{eqnarray}\label{eqB4}
&& \frac{1}{m} \sum_{j=1}^m X_{j,k}\bkg\xi_j^{T}\bkg\xi_j=
\frac{1}{m} \sum_{j=1}^m \ep[X_{j,k}\bkg\xi_j^{T}\bkg\xi_j|\Cal
F_{j-1}] +o(1) \non && = \frac{1}{m} \sum_{j=1}^m
\big(\ep[\pi_k(\bk x,\bkg\xi)\bkg\xi^{T}\bkg\xi]\big)\big|_{\bk
x=\widehat{\smallbkg\theta}_{j-1}} +o(1) \quad
a.s.
\end{eqnarray}
and
$$ \frac{1}{m}\sum_{j=1}^m X_{j,k}(Y_{j,k}-\bkg\theta_k\bkg\xi_j^{T})\bkg\xi_j\to 0 \quad a.s. $$
It follows that
\begin{align*}
\big\|\widehat{\bkg\theta}_{m,k}-\bkg\theta_k\big\|^2 & \left\{\min_{\|\bk y\|=1, \bk x\in \smallbkg\Theta_k}\bk y
\big(\ep[\pi_k(\bk x,\bkg\xi)\bkg\xi^{T}\bkg\xi]\big)\bk y^T-o(1)\right\}\\
\le &\big(\widehat{\bkg\theta}_{m,k}-\bkg\theta_k\big)\frac{1}{m}
[\sum_{j=1}^m X_{j,k}\bkg\xi_j^{T}\bkg\xi_j]\big(\widehat{\bkg\theta}_{m,k}-\bkg\theta_k\big)^T\\
=&\frac{1}{m}\sum_{j=1}^m
X_{j,k}(Y_{j,k}-\bkg\theta_k\bkg\xi_j^{T})\bkg\xi_j\big(\widehat{\bkg\theta}_{m,k}-\bkg\theta_k\big)^T
=o(1)\big\|\widehat{\bkg\theta}_{m,k}-\bkg\theta_k\big\| \quad a.s.
\end{align*}
Hence
\begin{equation}\label{eqB5}
\widehat{\bkg\theta}_{m,k}\to \bkg\theta_k\quad a.s.
\end{equation}
Now, by (\ref{eqB4}) and (\ref{eqB5}),
$$\frac{1}{m} \sum_{j=1}^m X_{j,k}\bkg\xi_j^{T}\bkg\xi_j=\bk I_{\xi k}+o(1) \quad
a.s., $$
which, together with (\ref{eqB3}), implies that
$$\widehat{\bkg\theta}_{m,k}-\bkg\theta_k
=\frac{1}{m}\sum_{j=1}^m
X_{j,k}(Y_{j,k}-\bkg\theta_k\bkg\xi_j^{T})\bkg\xi_j \bk
I_{\xi,k}^{-1}\big(1+o(1)\big) \;\; a.s.
$$
Notice $ \ep[(Y_{j,k}-
\bkg\theta_k\bkg\xi_j^{T}\big)\bkg\xi_j\big|\bkg\xi_j ]=0$.
(\ref{rateconsi}) is satisfied with $\bk
h_k=(Y_k-\bkg\theta_k\bkg\xi^{T})\bkg\xi \bk I_{\xi,k}^{-1}$.

\bigskip

{\bf Proof of Remark \ref{rek3.0}:} Now, the model is
$$ Y_j=X_{j,1}\mu_1+X_{j,2}\mu_2+\bkg\beta\widetilde{\bkg\xi}_j^T+\epsilon_j. $$
Let
$$ S(\mu_1^{\ast}, \mu_2^{\ast}, \bkg\beta^{\ast})=\sum_{j=1}^m
(Y_j-X_{j,1}\mu_1^{\ast}-X_{j,2}\mu_2^{\ast}-\bkg\beta^{\ast}\widetilde{\bkg\xi}_j^T)^2. $$
Wirte  $\bkg\eta_j=(X_{j,1},X_{j,2}, \widetilde{\bkg\xi})$ and $\bkg\theta=(\mu_1,\mu_2,\bkg\beta)$.
 With the same argument of proving Corollary \ref{rek3.1}, we have that the LSE $\widehat{\bkg\theta}_m$
 exists and is unique, and further,
satisfies the equation:
$$ (\widehat{\bkg\theta}_m-\bkg\theta)\left[\frac{1}{m}\sum_{j=1}^m \eta_j^T\eta_j\right]
=\frac{1}{m}\sum_{j=1}^m \epsilon_j\eta_j. $$
Write
$$ \bk I(\bk\theta)=\begin{pmatrix} \pi_1(\frac{\mu_1-\mu_2}{T})      & 0                       & \pi_1(\frac{\mu_1-\mu_2}{T}) \bk a \\
                                   0                   &\pi_2(\frac{\mu_1-\mu_2}{T})        & \pi_2(\frac{\mu_1-\mu_2}{T}) \bk a \\
                              \pi_1(\frac{\mu_1-\mu_2}{T})\bk a^T &\pi_2(\frac{\mu_1-\mu_2}{T})\bk a^T & \ep\widetilde{\bkg\xi}^T\bkg\xi
                 \end{pmatrix}. $$
Notice the assignment probability at stage $j$ depends only on the
estimator $\widehat{\bkg\theta}_{j-1}$ and does depend on the covariate
$\bkg\xi_j$. It follows that $v_1=\pi_1(\frac{\mu_1-\mu_2}{T})$, $v_2=\pi_2(\frac{\mu_1-\mu_2}{T})$, $\ep[\epsilon_j \bkg\eta_j\big|\Cal F_{j-1}]=0$,
$$\Var[\epsilon_j \bkg\eta_j\big|\Cal F_{j-1}]
=\sigma^2\ep[  \bkg\eta_j^T\bkg\eta_j\big|\Cal F_{j-1}]=\sigma^2\bk I(\widehat{\bk\theta}_{j-1})$$
and
$$\frac{1}{m}\sum_{j=1}^m \eta_j^T\eta_j=\frac{1}{m}\sum_{j=1}^m\bk I(\widehat{\bk\theta}_{j-1})+o(1). $$
So, similar to (\ref{eqB5}),
$$ \widehat{\bkg\theta}_m\to \bkg\theta\quad a.s. $$
Further
\begin{equation} \label{eqB6} \widehat{\bkg\theta}_m-\bkg\theta
=\frac{1}{m}\sum_{j=1}^m \epsilon_j\eta_j\bk I^{-1}\big(1+o(1)\big)\quad a.s.,
\end{equation}
where $\bk I=\bk I(\bkg\theta)$ and
$$\bk I^{-1}=\begin{pmatrix} 1/v_1+\bk a\bk
I_{\widetilde{\xi}}^{-1}\bk a^T & \bk a\bk I_{\widetilde{\xi}}^{-1}\bk a^T  & -\bk a I_{\widetilde{\xi}}^{-1}\\
\;\; \bk a\bk I_{\widetilde{\xi}}^{-1}\bk a^T  & 1/v_2+\bk a\bk
I_{\widetilde{\xi}}^{-1}\bk a^T & -\bk a I_{\widetilde{\xi}}^{-1} \\
-I_{\widetilde{\xi}}^{-1}\bk a^T & - \bk I_{\widetilde{\xi}}^{-1}\bk a^T & \bk I_{\widetilde{\xi}}^{-1}
\end{pmatrix}.$$
It follows that
$$ \sqrt{m}(\widehat{\bkg\theta}_m-\bkg\theta)
=\frac{1}{\sqrt{m}}\sum_{j=1}^m \epsilon_j\eta_j\bk I^{-1}+o_P(1)\overset{D}\to N(\bk 0,\bk I^{-1}). $$
From (\ref{eqB6}), it follows that
$$ \widehat{\mu}_{m1}-\widehat{\mu}_{m2}=\left(\frac{1}{n v_1}\sum_{j=1}^m \epsilon_jX_{j,1}
-\frac{1}{n v_2}\sum_{j=1}^m \epsilon_jX_{j,2}\right)(1+o(1))\quad a.s.$$
The remainder of the proof is similar to that of Theorem \ref{th2}.
\end{document}